\documentclass{amsart}
\usepackage{graphicx}
\theoremstyle{definition}

\theoremstyle{remark}

\numberwithin{equation}{section}

\begin{document}

\title{One-parameter Generalizations of Ramanujan's Formula for $\pi$}

\author {\textbf{Moa Apagodu}}
\address{ Department of Mathematics, Virginia Commonwealth
University, Richmond, VA 23284, USA}

\begin{abstract}
Several terminating generalizations of Ramanujan's formula for
$\frac{1}{\pi}$ with complete WZ proofs are given.

\end{abstract}


\maketitle

\noindent One of Ramanujan's [2] infinite series representation for
 $\frac{1}{\pi}$ is the series

$$
\sum_{k=0}^{\infty}
(-1)^k(4k+1)\frac{(\frac{1}{2})_k^3}{k!^3}=\frac{2}{\pi} \quad .
\eqno(Ramanujan)
$$

\noindent Zeilbeger [5] gave a short WZ proof of ({\it {Ramanujan}})
by first proving a one-parameter generalization, namely

$$
\sum_{k=0}^{\infty}(-1)^k(4k+1)\frac{(\frac{1}{2})_k^2(-n)_k}{k!^2(\frac{3}{2}+n)_k}=
\frac{\Gamma(\frac{3}{2}+n)}{\Gamma(\frac{3}{2})\Gamma(n+1)}
\,\,\,\quad , \eqno(Zeilberger)
$$

\vskip .2in

\noindent of ({\it {Ramanujan}}) for {\it nonnegative} integers $n$
using WZ method, and then evaluating both sides of the identity at
$n=-\frac{1}{2}$, thanks to Carlson's theorem (see below). \\

\noindent In this article, following Zeilberger's approach, we
provide several more one-parameter generalizations of ({\it
{Ramanujan}}) complete with their WZ proofs. These generalizations
(identities) are of interest on their own right as they appear to
be new at least for us. But first, \\

\noindent \textbf{Notation}: We denote a hypergeometric series\\

$$
{}_3 F_2 \left ( {{a,b,c} \atop {d,e}}  ; z, p(k) \right)
=\sum_{k=0}^{\infty}p(k)\frac{(a)_k(b)_k(c)_k}{k!(d)_k(e)_k}z^k
\quad ,
$$

\vskip .2in

\noindent by $F(a,b,c;d,e;z,p(k)),$ where
$(a)_j=\frac{\Gamma{(a+j)}}{\Gamma{(a)}}$ and $p$ is a polynomial in $k$.\\

\noindent Observe that the above series can also be viewed as a
${}_4F_3$ hypergeometric series. The following well-known theorem
due to Carlson is used to justify that if an identity holds for {\it
positive} integers,
then it also holds for rational arguments under suitable conditions. \\

\noindent \textbf{Theorem :} (Carlson [2]) If $f(z)$ is analytic and
is $\mathrm
{0}\left(e^{k\mathopen{\vert}z\mathclose{\vert}}\right)$, where
$k<\pi$, for $Re(z)\geq 0$, and if $f(z)=0$ for $z=0,1,2,\ldots$,
then $f(z)$ is identically zero.\\

\noindent \textbf{Theorem 1 :}
$$
\sum_{k=0}^{\infty}(-1)^k(4k+1){{(-n)_k^2(\frac{1}{2})_k}\over
{k!\left(n+\frac{3}{2}\right)_k^2}} =
\left(\frac{1}{4}\right)^n\frac{\left(\frac{3}{2}\right)_n^2}{\left(\frac{5}{4}\right)_n\left(\frac{3}{4}\right)_n}
\quad.
$$

\noindent \textbf{Proof :} \\

\noindent Let $F(n,k)$ be the summand divided by the right hand side
of the equality. Construct $G(n,k)=R(n,k)F(n,k),$ where $R(n,k)$ is
the rational function (certificate)

$$
R(n,k)=-{{(6n^2+10n+4+k-2k^2)k}\over {(n-k+1)^2(4k+1)}} \quad ,
$$

\vskip .1in

\noindent so that $F(n+1,k)-F(n,k)=G(n,k+1)-G(n,k)$. Now sum both
sides of this last equation with respect to k ($k=0$ to $k=\infty$),
to see that the right hand side telescopes to zero from which it
follows that $\sum F(n,k) = Constant$. Finally, plugging in $n=0$,
we get $\sum_{k=0}^{\infty}F(n,k)=1$ completing the WZ proof of the
theorem for {\it nonnegative} integers $n$. To deduce ({\it
Ramanujan}), substitute $n=-\frac{1}{2}$ which is legitimate by
Carlson's theorem. QED

\vskip .2in

\noindent In our notation, the statement of theorem 1 is equivalent
to

$$
F(-n,-n,\frac{1}{2};1,n+\frac{3}{2};-1,4k+1)=
\left(\frac{1}{4}\right)^n\frac{\left(\frac{3}{2}\right)_n^2}{\left(\frac{5}{4}\right)_n\left(\frac{3}{4}\right)_n}
\quad.
$$

\vskip .2in

\noindent Below we provide more terminating generalizations which
reduces to ({\it Ramanujan}) when evaluated at $n=-{{1}\over {2a}}$,
where $a$ is the coefficient of $n$ in $F(-an,b,c;d,e;z,p(k))$. In
the remaining generalizations except theorem 2, the right hand side
do not automatically simplify to $\frac{2}{\pi}$ which by itself
gives some interesting relationship between different Gamma and
trigonometric values. To wit, in theorem 6 below, when we evaluate
the right hand side of the identity at $n=-\frac{1}{2}$, we get
$\frac{\sqrt{5}}{\pi (\cos(\frac{\pi}{5})+ \cos(\frac{2\pi}{5}))}$
which equals $\frac{2}{\pi}$ (To see $\cos(\frac{\pi}{5})+
\cos(\frac{2\pi}{5})=\frac{1}{2}\sqrt{5}$, consider roots of $4x^2-2x-1=0)$.\\

\noindent \textbf{Theorem 2 :}
$$
F\left(-n,-2n-\frac{1}{2};\frac{1}{2},n+\frac{3}{2};2n+2,-1,4k+1\right)=
\left(\frac{2^2}{3^3}\right)^n\frac{\left(\frac{3}{2}\right)_n^2}{\left(\frac{4}{3}\right)_n\left(\frac{2}{3}\right)_n}
\quad .
$$

\noindent \textbf{Proof :} Let\\

\noindent
$R(n,k)=(184n^4+658n^3-44k^2n^2+22kn^2+868n^2+38kn-76k^2n+500n+106+4k^4+17k-4k^3-33k^2)\times$ \\

$$
\frac{2k}{(4k+1)(4n+5-2k)(4n+3-2k)(-k+n+1)(2n+2+k)} \quad ,
$$

\noindent and proceed as in Theorem 1.

 \vskip .2in

\noindent \textbf{Theorem 3 :}

$$
F\left(-2n,-n+\frac{1}{4},\frac{1}{2};n+\frac{5}{4},2n+\frac{3}{2};-1,4k+1\right)=\left(\frac{2^2}{3^3}\right)^n
\frac{\left(\frac{5}{4}\right)_n^2}{\left(\frac{13}{12}\right)_n\left(\frac{5}{12}\right)_n}\quad.
$$

\noindent \textbf{Proof :} Let \\

\noindent
$R(n,k)=(2944n^4+7584n^3+352kn^2-704k^2n^2+7096n^2+432kn-864k^2n+2846n+64k^4+142k-64k^3-268k^2+411)\times$

$$
\frac{-k}{4\left(4k+1\right)\left(2n-k+2\right)\left(2n-k+1\right)\left(4n+3-4k\right)\left(4n+3+2k\right)}
\quad ,
$$

\noindent and proceed as in Theorem 1.

\vskip .2in

\noindent \textbf{Theorem 4:}
$$
F\left(-n,-3n-1,\frac{1}{2};n+\frac{3}{2},3n+\frac{5}{2};-1,4k+1\right)=\left(\frac{3^3}{2^8}\right)^n
\frac{\left(\frac{7}{6}\right)_n\left(\frac{5}{6}\right)_n\left(\frac{3}{2}\right)_n^2}{\left(\frac{5}{8}\right)_n\left(\frac{7}{8}\right)_n\left(\frac{9}{8}\right)_n\left(\frac{11}{8}\right)_n}
\quad .
$$

\noindent \textbf{Proof :} Let \\

\noindent
$R(n,k):=(8470+58774n+963k+210k^4-8460k^2n+7138kn^2+4286kn-1872k^2+251126n^3+12k^5-215k^3+208908n^4-8k^6-10528k^2n^3+1452kn^4+5264kn^3-
14214k^2n^2+91488n^5-448k^3n+167522n^2-2904k^2n^4-248k^3n^2+16488n^6+248k^4n^2+448k^4n)\times$\\

$$
\frac{-k}{\left(4k+1\right)\left(3n+4-k\right)\left(3n+3-k\right)\left(3n+2-k\right)\left(-k+n+1\right)\left(6n+7+2k\right)\left(6n+5+2k\right)}
$$

\noindent and proceed as in Theorem 1.

\vskip .2in

\noindent \textbf{Theorem 5 :}

$$
F\left(-3n,-n+\frac{1}{3},\frac{1}{2};n+\frac{7}{6},3n+\frac{3}{2};-1,4k+1\right)=\left(\frac{3^3}{2^8}\right)^n{
{\left(\frac{5}{6}\right)_n\left(\frac{1}{2}\right)_n\left(\frac{7}{6}\right)_n^2}\over
{\left(\frac{25}{24}\right)_n\left(\frac{7}{24}\right)_n\left(\frac{13}{24}\right)_n\left(\frac{19}{24}\right)_n}}\quad.
$$

\noindent \textbf{Proof :} Let\\

\noindent
$R(n,k)=(33736+412476n+12183k+7146k^4-167112k^2n+230202kn^2+86418kn-22458k^2+5023782n^3+972k^5-
7551k^3+6796548n^4-648k^6-539136k^2n^3+117612kn^4+269568kn^3-455382k^2n^2+4739472n^5-22896k^3n+
2011788n^2-235224k^2n^4-20088k^3n^2+1335528n^6+20088k^4n^2+22896k^4n)\times$

$$
\frac{-k/\left(27\left(4k+1\right)\right)}{\left(3n+2-3k\right)\left(3n+3-k\right)\left(3n+2-k\right)\left(3n-k+1\right)\left(6n+5+2k\right)\left(6n+3+2k\right)}
$$

\noindent and proceed as in Theorem 1.

\vskip .2in

\noindent \textbf{Theorem 6 :}

$$
F\left(-n,-4n-\frac{3}{2},\frac{1}{2};n+\frac{3}{2},4n+3;-1,4k+1\right)=
\left(\frac{2^8}{5^5}\right)^n\frac{\left(\frac{5}{4}\right)_n\left(\frac{3}{4}\right)_n\left(\frac{3}{2}\right)_n^2}
{\left(\frac{6}{5}\right)_n\left(\frac{7}{5}\right)_n\left(\frac{3}{5}\right)_n\left(\frac{4}{5}\right)_n}
\quad .
$$

\noindent \textbf{Proof :} Let\\

\noindent
$R(n,k):=(1360680+12454594n+118956k+22445k^4-1561554k^2n+2189040kn^2+792996kn-232109k^2+
110282888n^3+1528k^5+22139008n^7-23084k^3+152482100n^4+2937856n^8-2048k^6n+1656k^5n^2+3072k^5n+
16k^8-1000k^6-6375798k^2n^3-1104k^6n^2+1133632kn^5+203392kn^6+2617070kn^4-406784k^2n^6+3201812kn^3-
4339096k^2n^2+133488542n^5-97496k^3n+49304668n^2-2267264k^2n^5-5226636k^2n^4-32k^7-
111304k^3n^3-155798k^3n^2+72279728n^6-30016k^3n^4+111304k^4n^3+155108k^4n^2+96216k^4n+30016k^4n^4))\times
$

$$
\frac{-2k/((4k+1)(8n+11-2k)(8n+9-2k))}{(8n+7-2k)(8n+5-2k)(-k+n+1)(4n+5+k)(4n+k+4)(4n+k+3)}
$$

\noindent and proceed as in Theorem 1.

\vskip .2in

\noindent \textbf{Theorem 7 :}
$$
F\left(-3n,-2n+\frac{1}{6},\frac{1}{2};2n+\frac{4}{3},3n+\frac{3}{2};-1,4k+1\right)=
\left(\frac{2^23^3}{5^5}\right)^n\frac{\left(\frac{5}{6}\right)_n\left(\frac{1}{2}\right)_n\left(\frac{7}{6}\right)_n^2}{\left(\frac{7}{15}\right)_n\left(\frac{13}{15}\right)_n\left(\frac{16}{15}\right)_n\left(\frac{4}{15}\right)_n}
\quad .
$$

\noindent \textbf{Proof :} Let \\

\noindent
$R(n,k)=(4701560+76180392n+1024494k+708705k^4-22123845k^2n+56137239kn^2+11776977kn-1840161k^2+2005650450n^3+
248832k^5+2986094808n^7-814086k^3+4671194832n^4+633425184n^8-520992k^6n+717336k^5n^2+781488k^5n+11664k^8-
152280k^6-277628958k^2n^3-478224k^6n^2+144102888kn^5+43337592kn^6+196729884kn^4-86675184k^2n^6+141103782kn^3-
108466425k^2n^2+6791227920n^5-5655312k^3n+524305530n^2-288205776k^2n^5-391357332k^2n^4-23328k^7-18314424k^3*n^3-
15172434k^3n^2+6025575744n^6-8409744k^3n^4+18314424k^4n^3+14873544k^4n^2+5329692k^4n+8409744k^4n^4)\times$

$$
\frac{-2k/(27(4k+1)(12n+11-6k)(12n+5-6k))}{(3n+3-k)(3n+2-k)(3n-k+1)(6n+5+2k)(6n+3+2k)(6n+4+3k)}
$$

\noindent and proceed as in Theorem 1.

\vskip .2in

 \noindent \textbf{Theorem 8 :}
$
F\left(-4n,-n+\frac{3}{8},\frac{1}{2};n+\frac{9}{8},4n+\frac{3}{2};-1,4k+1\right)=\left(\frac{2^8}{5^5}\right)^n
\frac{\left(\frac{7}{8}\right)_n\left(\frac{3}{8}\right)_n\left(\frac{9}{8}\right)_n^2}
{\left(\frac{33}{40}\right)_n\left(\frac{41}{40}\right)_n\left(\frac{9}{40}\right)_n\left(\frac{17}{40}\right)_n}
\quad . $ \\

\noindent \textbf{Proof :} Let \\

\noindent $R(n,k)=
(12815055+232274998n+4590732k-396544k^6-1249280k^6n-1130496k^6n^2+28224057344n^5+7320264608n^3+
18289568000n^4+26350223360n^6+54254040kn+208273408kn^6+942632960kn^4+3008364544n^8+671602688kn^3+
692224000kn^5+263857056kn^2-32768k^7-103023312k^2n-8360336k^2-416546816k^2n^6-1877581824k^2n^4-1326237696k^2n^3-1384448000k^2n^5-
513366592k^2n^2+21002112k^4n+2969696k^4+30736384k^4n^4+67870720k^4n^3+56542208k^4n^2-21782912k^3n-3231872k^3-30736384k^3n^4-
67870720k^3n^3-57248768k^3n^2+16384k^8+1761550336n^2+13645250560n^7+1873920k^5n+623488k^5+1695744k^5n^2)\times $\\

$$
\frac{-k/{(128(4k+1)(8n+5-8k)(4n+4-k))}}{(4n+3-k)(4n-k+2)(4n-k+1)(8n+7+2k)(8n+5+2k)(8n+3+2k)}
$$

\vskip .2in

\noindent \textbf{Theorem 9 :}

$$
F\left(-2n,-3n-\frac{1}{4},\frac{1}{2};2n+\frac{3}{2},3n+\frac{7}{4};-1,4k+1\right)=
\left(\frac{2^23^3}{5^5}\right)^n\frac{\left(\frac{11}{12}\right)_n\left(\frac{7}{12}\right)_n\left(\frac{5}{4}\right)_n^2}{\left(\frac{11}{20}\right)_n\left(\frac{19}{20}\right)_n\left(\frac{23}{20}\right)_n\left(\frac{7}{20}\right)_n}
\quad .
$$

\noindent \textbf{Proof :} Let \\

\noindent
$R(n,k)=(22623909+306149258n+3494086k+14411873792n^5+5772117536n^3+11505823872n^4+
232853504kn^5+367020544kn^4+305077760kn^3+34504208kn+60874752kn^6+141134368kn^2+
11083683840n^6+889749504n^8-6486428k^2-46570e008k^2n^5-731087872k^2n^4-602739712k^2n^3-
65967264k^2n-121749504k^2n^6-275188800k^2n^2-1968960k^3-11812864k^3n^4-29663232k^3n^3-
12059136k^3n-28235776k^3n^2+1779904k^4+29663232k^4n^3+11812864k^4n^4+11531776k^4n+
27815936k^4n^2+448000k^5+1265664k^5n+1007616k^5n^2+
1775873160n^2+4787625984n^7-279552k^6-843776k^6n-671744k^6n^2-32768k^7+16384k^8) \times$\\

$$
\frac{-k/(4(4k+1)(12n+13-4k)(12n+9-4k)(12n+5-4k))}{(2n-k+2)(2n-k+1)(4n+3+2k)(12n+11+4k)(12n+7+4k)}
$$

\noindent and proceed as in theorem 1.
 \vskip .2in

\noindent A similar proof can be constructed for the following two
identities using Zeilberger algorithm. \\

\noindent \textbf{Theorem 10 :}

$$
F\left(-4n,-3n+\frac{1}{8},\frac{1}{2};3n+\frac{11}{8},4n+\frac{3}{2};-1,4k+1\right)=\left(\frac{2^83^3}{7^7}\right)^n
\frac{\left(\frac{11}{24}\right)_n\left(\frac{3}{8}\right)_n\left(\frac{7}{8}\right)_n\left(\frac{19}{24}\right)_n\left(\frac{9}{8}\right)_n^2}
{\left(\frac{11}{56}\right)_n\left(\frac{43}{56}\right)_n\left(\frac{19}{56}\right)_n\left(\frac{51}{56}\right)_n\left(\frac{27}{56}\right)_n\left(\frac{59}{56}\right)_n}
\quad .
$$

\vskip .2in

\noindent \textbf{Theorem 11 :}

$$
F\left(-3n,-4n-\frac{1}{6},\frac{1}{2};3n+\frac{3}{2},4n+\frac{5}{3};-1,4k+1\right)=\left(\frac{2^83^3}{7^7}\right)^n
\frac{\left(\frac{11}{12}\right)_n\left(\frac{5}{6}\right)_n\left(\frac{1}{2}\right)_n\left(\frac{5}{12}\right)_n\left(\frac{7}{6}\right)_n^2}
{\left(\frac{11}{21}\right)_n\left(\frac{23}{21}\right)_n\left(\frac{5}{21}\right)_n\left(\frac{17}{21}\right)_n\left(\frac{8}{21}\right)_n\left(\frac{20}{21}\right)_n}
\quad .
$$

\vskip .2in

\noindent \textbf{Conclusion}\\

\noindent In this article we considered one parameter
generalizations of one of the many formulas of Ramanujan for $\pi$.
It would be interesting to find if similar generalizations exist for
other similar formulas for $\frac{1}{\pi}$. For example a notable one is the series\\

$$
2\sqrt{2}\sum_{k=0}^{\infty}
\left(\frac{1}{99}\right)^{4k+2}(1103+26390k)\frac{(\frac{1}{4})_k
(\frac{1}{2})_k (\frac{3}{4})_k}{k!^3}=\frac{1}{\pi}  \quad .
$$

\noindent See http://mathworld.wolfram.com/PiFormulas.html for
complete list of similar formulas.\\

\noindent \textbf{Acknowledgement}\\

\noindent I would like to thank Doron Zeilberger for introducing me
to the intriguing question posed by John Greene (to find a
one-parameter family generalizing the exact evaluation, in terms of
the Gamma function evaluated at rational arguments, of $F({{1} \over
{4}}, {{3} \over {4}},1, - {{1} \over {63}})$ ).\\

\noindent {\bf References}\\

\bibliographystyle{amsplain}

\noindent [1] G.E. Andrews, R. Askey, and R. Roy, {\it ``Special
Functions''}, Cambridge Univ. Press, 1999.\\

\noindent [2] George E. Andrews, Richard A. Askey, Bruce C. Berndt,
K.G. Ramanathan, and Robert A. Rankin, {\it Ramanujan Revisited},
Proceeding of The Centenary Conference, University of Illinois at
Urbana-Champaign, June 1-5, Academic Press (1987).\\

\noindent [3] W.N. Bailey, {\it ``Generalized Hypergeometric
Series''}, Cambridge University Press, 1935. Reprinted by Hafner
Pub. Co., New York, 1972.\\

\noindent [4] I. Gessel and D. Stanton, {\it Strange evaluations of
hypergeometric series}, SIAM J. Math. Anal. {\bf13}(1982),
295-308.\\

\noindent [5] M. Mohammed and D. Zeilberger, {\it Sharp upper bounds
for the orders outputted by the Zeilberger and q-Zeilberger
algorithms}, J. Symbolic Computation {\bf 39} (2005), 201-207.\\

\noindent [6] S. B. Ekhad and D. Zeilberger, {\it A WZ proof of
Ramanujan's formula for $\pi$}, Geometry, Analysis, and
Mechanics, ed. by J.M. Rassias, World Scientific, Singapore, 1994, 107-108.\\

\end{document}